\documentclass[a4paper, conference]{IEEEtran}
\usepackage{amsmath}
\usepackage{amssymb}
\usepackage{amsfonts}
\usepackage{graphicx}
\usepackage{subcaption}
\usepackage{epstopdf}
\usepackage{pseudocode}
\usepackage{algorithm}
\usepackage{enumerate}
\usepackage{url}

\def\lb{\lambda} 
\def\bu{{\boldsymbol u}} 
\def\bx{{\boldsymbol x}}
\def\by{{\boldsymbol y}}
\def\ba{{\boldsymbol a}}
\def\bb{{\boldsymbol b}}

\def\b0{{\boldsymbol 0}}

\def\bp{{\boldsymbol p}}
\def\bpsi{{\boldsymbol \psi}}
\def\bdelta{{\boldsymbol \delta}}
\def\CN{{\cal N}}

\def\eps{\varepsilon}

\def\be{\begin{equation}}
\def\ee{\end{equation}}
\def\bea{\begin{eqnarray}}
\def\eea{\end{eqnarray}}

\def\nn{{\nonumber}}

\title{Hamiltonian System Approach to \\ Distributed Spectral Decomposition in Networks
}

\author{\IEEEauthorblockN{Konstantin Avrachenkov}
\IEEEauthorblockA{INRIA, France\\
Email: k.avrachenkov@inria.fr}
\and
\IEEEauthorblockN{Philippe Jacquet}
\IEEEauthorblockA{Nokia Bell Labs, France \\
Email: philippe.jacquet@nokia-bell-labs.com}
\and
\IEEEauthorblockN{Jithin K.\ Sreedharan}
\IEEEauthorblockA{Purdue University, USA\\
Email: jithinks@purdue.edu}}

\begin{document}

\maketitle
\thispagestyle{empty}
\pagestyle{empty}

\begin{abstract}
Because of the significant increase in size and complexity of the networks, the distributed computation of eigenvalues and eigenvectors of graph matrices has become very challenging and yet it remains as important as before. In this paper we develop efficient distributed algorithms to detect, with higher resolution, closely situated eigenvalues and corresponding eigenvectors of symmetric graph matrices. We model the system of graph spectral computation as physical systems with Lagrangian and Hamiltonian dynamics. The spectrum of Laplacian matrix, in particular, is framed as a classical spring-mass system with Lagrangian dynamics. The spectrum of any general symmetric graph matrix turns out to have a simple connection with quantum systems and it can be thus formulated as a solution to a Schr\"odinger-type differential equation. Taking into account the higher resolution requirement in the spectrum computation and the related stability issues in the numerical solution of the underlying differential equation, we propose the application of symplectic integrators to the calculation of eigenspectrum. The effectiveness of the proposed techniques is demonstrated with numerical simulations on real-world networks of different sizes and complexities.

\end{abstract}

\section{Introduction}
Consider an undirected graph $G=(V,E)$ with $V$ as the vertex set ($n:= \lvert V \rvert$) and $E$ as the edge set ($m:= \lvert E \rvert$). Let $M$ be any symmetric matrix associated with $G$. Broadly speaking, we say that
$M$ is a graph matrix if it has non-zero elements on the edge set. Due to symmetry, the eigenvalues of $M$ are real and can be ranked in ascending order as $\lb_1 \le \lb_2\le \ldots \le\lb_n$. We investigate efficient techniques to compute the eigenvalues $\lambda_1,\ldots, \lambda_n$ and the corresponding eigenvectors $\bu_1,\ldots,\bu_n$, in a distributed way.

We define two typical matrices which appear frequently in network analysis. First one is the adjacency matrix $A \in \mathbb{R}^{n \times n}$ in which the individual entries are given by
\[
a_{uv} =
\left\{ \begin{array}{ll}
1 & \mbox{if} \ u \ \mbox{is a neighbour of} \ v,\\
0 & \mbox{otherwise}.
\end{array}\right.
\]
Since the focus in this paper is on undirected graphs, $A^{\intercal}=A$. The matrix $A$ is also called the unweighted adjacency matrix and one can also define a weighted version in which the weight $1$ for an edge is replaced by any weight such that $a_{uv} = a_{vu}$. Another matrix which is found very common in many graph theoretic problems is the Laplacian matrix $L= [\ell_{i,j}]:= D-A$. Here the matrix $D$ is a diagonal matrix with diagonal elements equal to degrees of the nodes $d_1, \ldots, d_n$.

\subsection{Applications of graph spectrum}
The knowledge of $\{\lambda_i\}$'s and $\{\bu_i\}$'s can be made use of in many ways. For instance,
spectral clustering is a prominent solution which exploits the first $k$ eigenvectors of the Laplacian matrix for identifying the clusters in a network \cite{VonLuxburg07}.
Another classical use of Laplacian eigenvalues is in computing the number of spanning trees of a graph $G$ which is ${n^{-1}\,\lb_2 \lb_3 \ldots \lb_n}$. Many studies have been conducted on the use of the spectrum of the adjacency matrix such as computing the number of traingles of a netwok (both locally and globally) \cite{Tsourakakis_ICDM08}, graph dimensionality reduction and link reduction \cite{Kempe04} etc.
%

Two applications relevant to multi-agent and multi-dimensional systems will be explained in detail in Section~\ref{sec:applications}.

\subsection{Our basic approach}
%
Let $M$ be any symmetric graph matrix. Consider the following Schr\"odinger-type differential equation
\be
\frac{\partial}{\partial t}\bpsi(t) = iM\,\bpsi(t),
\label{eq:de_infocom16}
\ee
where $\bpsi(t)$ is a complex valued $n$ dimension vector, which can be interpreted as the wave function of a hypothetical quantum system. The solution of this differential equation with the boundary condition $\bpsi(0)=\ba_0$ is $\exp (iMt) \ba_0$. Subjecting this solution to the Fourier transform provides a decomposition in terms of eigenvalues and eigenvectors as follows:
\bea
\int_{-\infty}^{+\infty}e^{i M t}\ba_0 e^{-it\theta}dt = 2 \pi \sum_{j=1}^n \delta_{\lambda_j}(\theta)\bu_j (\bu_j^\intercal \ba_0),
\label{eq:dirac_decompsn}
\eea
where $\delta_{\lambda_j}$ is the Dirac function shifted by $\lambda_j$. This follows from the eigen-decomposition of the matrix $M$, $e^{iMt} = \sum_j e^{i t \lambda_j} \bu_j \bu_j^{\intercal}$. In order to avoid the harmonic oscillations created from finite and discretized version of the above Fourier transform, which will mask the Dirac peaks, the following Gaussian smoothing can be performed. For $v>0$,
\bea
\lefteqn{\frac{1}{2\pi} \int_{-\infty}^{+\infty}e^{i M t} \ba_0 e^{-t^2v/2}e^{-it\theta}dt}\nn \\
&=&\sum_{j=1}^n \frac{1}{\sqrt{{2 \pi} v}} \exp(-\frac{(\lambda_j-\theta)^2}{2v})\bu_j (\bu_j^\intercal \ba_0).
\label{eq:gauss_approxn}
\eea
The right-hand side of the above expression leads us to a plot at each node $k$ with Gaussian peaks at each of the eigenvalues and the amplitude of the peak at $j$th eigenvalue is ${(\sqrt{2 \pi v})}^{-1} (\bu_j^{\intercal} \ba_0)\bu_j(k)$, proportional to the $k$th component in the eigenvector $\bu_j$.

The idea is to estimate the solution of \eqref{eq:de_infocom16} at $\eps$ intervals of time for a total of $s$ samples. One way to form an estimate to the left-hand side to \eqref{eq:gauss_approxn} is:
\begin{equation}
\eps \Re \Big( \bb_0+2\sum_{\ell= 1}^{s} e^{i\eps \ell M} {\bb}_0 e^{-i \ell \eps \theta}e^{-\ell^2\eps^2v/2} \Big).
\label{eq:f_expn}
\end{equation}

In \cite{AvJaSr_Infocom16}, we use the above approximation with various approaches based on diffusion, gossiping and quantum random walks for distributed computation of the eigenvalues and the eigenvectors. Let us discuss some challenges in this approach.

\subsection*{Issues in the computation}
While doing numerical experiments, we have observed that the approaches in \cite{AvJaSr_Infocom16} work well for larger eigenvalues of the adjacency matrix of a graph but they do not perform that well when one needs to distinguish between the eigenvalues which are very close to each other. One of the main techniques proposed there to solve \eqref{eq:de_infocom16} and to find  $e^{i\eps \ell M}$ in \eqref{eq:f_expn} are via $r$th order Runge-Kutta method and its implementation as a diffusion process in the network. The $r$-th order Runge-Kutta method has the convergence rate of $\mathcal{O}(\eps^r)$. We have observed that this is the case while checking the trajectory of the associated differential equation solution; the solution diverges, and it happens when a large number of iterations $s$ is required (see Section~\ref{sec:numerical_res}).


A larger value for $s$ is anticipated from our approximation in \eqref{eq:f_expn} due to the following facts. From the theory of Fourier transform and Nyquist sampling, the following conditions must be satisfied:
\be
\eps \leq \frac{\pi}{\lambda_{n}} \text{ and } s \geq \frac{2 \pi}{\eps \lambda_{\text{diff}}},
\label{eq:ft_constraints}
\ee
where $\lambda_{\text{diff}}$ is the maximum resolution we require in the frequency (eigenvalue) domain, which is ideally $\min_{i} |\lambda_i -\lambda_{i+1}|$. This explains that when dealing with graph matrices with larger $\lambda_n$ and require higher resolution, $s$ will take large values. For instance, in case of the Laplacian matrix, where the maximum eigenvalue is bounded as $\frac{n}{n-1} \Delta(G) \leq \lambda_{n} \leq 2 \Delta(G),$ with $\Delta(G)$ as the maximum degree of the graph and the lower eigenvalues are very close to each other, $s$ turns out to be typically a very large value.

In what follows, we propose solutions to the above mentioned issues.

\subsection{Related works}
The general idea of using mechanical oscillatory behaviour for the detection of eigenvalues has appeared in a few previous works, see e.g., \cite{franceschelli2013decentralized,sahai2012hearing}.
Though the technique in \cite{franceschelli2013decentralized} is close to ours, our methods differ by focussing on a Schr\"odinger-type equation and numerical integrators specific to it. Moreover, we demonstrate the efficiency and stability of the methods in real-world networks of varying sizes, in contrast to a small size synthetic network considered in \cite{franceschelli2013decentralized}, and our methods can be used to estimate eigenvectors as well.

In comparison to \cite{sahai2012hearing} we do not deform the system and we use new symplectic numerical integrators \cite{blanes2006symplectic,blanes2008splitting}. For the problem of distributed spectral decomposition in networks, one of the first and prominent works appeared in \cite{Kempe04}. But their algorithm requires distributed orthonormalization at each step and they solve this difficult operation via random walks. But if the graph is not well-connected (low conductance), this task will take a very long time to converge. Our distributed implementation based on fluid diffusion in the network does not require such orthonormalization.

\subsection{Contributions}
We make the following contributions and significantly improve the algorithms from our previous work:
\begin{enumerate}
\item We observe from our previous studies that the stability in trajectory of the differential equation solver is of significant influence in the eigenvalue-eigenvector technique based on \eqref{eq:de_infocom16}. Thus, we resort to geometric integrators to ensure the stability. In particular, by modeling the problem as a Hamiltonian system, we use symplectic integrators (SI) which protect the volume preservation of Hamiltonian dynamics, thus preserve stability and improve accuracy.
\item We propose algorithms that are easy to design without involving many parameters with interdependence, compared to the algorithms proposed in \cite{AvJaSr_Infocom16}.
\end{enumerate}

In the rest of the paper for clarity of presentation we mostly concentrate on the Laplacian matrix $L$ as an example for graph matrix $M$. We design algorithms based on Lagrangian as well as Hamiltonian dynamics, to compute the smallest $k$ eigenvalues and the respective eigenvectors of the Laplacian matrix efficiently. For simplicity, in this paper we do not consider Gaussian smoothing \eqref{eq:gauss_approxn}, but the proposed algorithms can be readily extended to include it.

The paper is organized as follows. In Section \ref{sec:mechanical_lag}, we explain a mass-spring analogy specific to Laplacian matrix and derive a method to identify the spectrum. Section \ref{sec:ham_dynamics} focuses on general symmetric matrices and develop techniques based on solving the Schr\"odinger-type equation efficiently. Section \ref{sec:distributed_impln} details a distributed implementation of the proposed algorithm. Section \ref{sec:numerical_res} contains numerical simulations on networks of different sizes. Section~\ref{sec:applications} contains two relevant applications to multi-agent and multi-dimensional systems. Section \ref{sec:conclusions} concludes the paper.

For convenience we summarize the important notation used in this paper in Table \ref{tab:list_symbols}.
\vspace*{-0.25cm}
\begin{table}[h]
\begin{center}
\begin{tabular}[h!]{|c||l|}\hline
  Notation & Meaning\\ \hline \hline
  $G$, $(V,E)$ & Graph, Node set and edge set \\ \hline
  $n$, $m$ & No.\ of nodes, no.\ of edges \\ \hline
  $A$ & Adjacency matrix \\ \hline
  $L$ & Laplacian matrix\\ \hline
  $\lambda_1, \ldots, \lambda_k$ & Smallest $k$ eigenvalues of $L$ in ascending order  \\ \hline
  $\bu_1, \ldots, \bu_k$ & Eigenvectors corresponding to $\lambda_1, \ldots, \lambda_k$ \\ \hline
  $d_j$ & Degree of node $j$ without including self loop \\ \hline
  $\Delta(G)$ & $\max\{d_1,\ldots,d_n\}$ \\ \hline
  $\CN(m)$ & Neighbor list of node $m$ without including self loop \\ \hline
  $\eps$ & Sampling interval in time domain \\ \hline
  $T,s$ & Total time frame and no.\ of samples $\lceil T/\eps \rceil$ \\ \hline
  $\bx(t), \bx_i$ & vector $\bx$ with index in continuous and discrete time \\ \hline
  $\bx_i[k]$ & $k$th component of a vector $\bx_i$ \\ \hline
\end{tabular}
\vspace*{0.2cm}
\caption{List of important notations}
\label{tab:list_symbols}
\end{center}
\vspace*{-0.75cm}
\end{table}

\section{Mechanical spring analogy with Lagrangian dynamics}
\label{sec:mechanical_lag}
Consider a hypothetical mechanical system representation of the graph $G$ in which unit masses are placed on the vertices and the edges are replaced with mechanical springs of unit stiffness. Using either Lagrangian or Netwonian mechanics, the dynamics of this system is described by the following system of differential equations
\be
\label{eq:LargDiff}
\ddot{\bx}(t) + L \bx(t) = \b0.
\ee
The system has the Hamiltonian function as ${\cal H} = \frac{1}{2} \dot{\bx}^\intercal I \dot{\bx} + \frac{1}{2} \bx^\intercal L \bx$, $I$ being the identity matrix of order $n$.

We note that once we obtain by some identification method the frequencies $\omega_k$ of the above oscillatory system, the eigenvalues of the Laplacian $L$ can be immediately retrieved by the simple formula $\lambda_k = |\omega_k^2|$. This will be made clearer later in this section.

Starting with a random initial vector $\bx(0)$, we can simulate the motion of this spring system. For the numerical integration, the Leapfrog or Verlet method \cite{verlet1967computer} technique can be applied.
Being an example of geometric integrator, Verlet method has several remarkable properties. It has the same computational complexity as the Euler method but it is of second order method (Euler method is employed in \cite{AvJaSr_Infocom16} as the first order distributed diffusion). In addition, the Verlet method is stable for oscillatory motion and conserves the errors in energy and computations \cite[Chapter~4]{leimkuhler2004simulating}. It has the following two forms. Let $\bp(t):= \dot{\bx}(t)$ and $\bx_i$ be the approximation of $\bx(i\eps)$, similarly $\bp_i$ for $\bp(i\eps)$. Here $\eps$ is the step size for integration. First, define
\[\bp_{1/2} = \bp_0 + \eps/2 (-L \bx_0).\]
Then, perform the following iterations
\bea
\bx_i & = & \bx_{i-1} + \eps \bp_{i-1/2} \nn \\
\bp_{i+1/2} & = & \bp_{i-1/2} + \eps (-L \bx_i). \nn
\eea
Equivalently, one can do the updates as
\bea
\bx_{i+1} & = & \bx_{i} + \eps p_{i} + \eps^2/2 (-L \bx_i) \nn \\
\bp_{i+1} & = & \bp_{i} + \eps [(-L \bx_i) + (-L \bx_{i+1})].\nn
\eea
We name the above algorithm as Order-2 Leapfrog.

Solution of the differential equation \eqref{eq:LargDiff} subject to the boundary values $\bx(0)=\ba_0$ and $\bp(0)=\bb_0$ is
\be
\bx(t) = \left(\frac{1}{2}\ba_0-i\frac{\bb_0}{\sqrt{\Lambda}} \right) e^{i t \sqrt{L}} +
 \left(\frac{1}{2}\ba_0+i\frac{\bb_0}{\sqrt{\Lambda}} \right) e^{-i t \sqrt{L}}, \nn
\ee
where we assume the decomposition of $L$ based on spectral theorem, i.e., $L=U \Lambda U^\intercal$ with $U$ as the orthonormal matrix with columns as eigenvectors and $\Lambda$ as the diagonal matrix formed from the eigenvalues. Further simplification of the above expression along the fact that $f(L)=U f(\Lambda) U^\intercal$, for any function $f$ which can be expressed in terms of power series, gives
\be
\bx(t) = \cos(t\sqrt{L}) \ba_0 + (\sqrt{L})^{-1} \sin (t \sqrt{L}) \bb_0. \nn
\ee
or $k$th component of $\bx(t)$ is
\[\ba_0[k] \cos (t\sqrt{\lambda_k})+\frac{\bb_0[k]}{\sqrt{\lambda_k}} \sin(t \sqrt{\lambda_k}).\]
Now we have
\bea
\lefteqn{\int_{-\infty}^{+\infty} \bx(t) e^{-it\theta}dt} \nn \\
& = & \int_{-\infty}^{+\infty} \sum_{k=1}^n \cos (t\sqrt{\lambda_k}) \bu_k (\bu_k^{\intercal} \ba_0) e^{-it\theta}dt \nn \\
& & +\int_{-\infty}^{+\infty} (\sqrt{L})^{-1} \sum_{k=1}^n \sin(t\sqrt{\lambda_k}) \bu_k (\bu_k^{\intercal} \bb_0) e^{-it\theta}dt \nn \\
& = & \sum_{k=1}^n \bu_k (\bu_k^{\intercal} \ba_0) \left(\pi[\delta(\theta-\sqrt{\lambda_k})+\delta(\theta+\sqrt{\lambda_k})] \right) \nn \\
& & \hspace*{- 0.7 cm} +(\sqrt{L})^{-1} \bu_k (\bu_k^{\intercal} \bb_0) \left(-\pi i [\delta(\theta-\sqrt{\lambda_k})-\delta(\theta+\sqrt{\lambda_k})] \right).\nn
\eea
Taking the real and positive spectrum will give $\pi \sum_{k=1}^n \bu_k (\bu_k^{\intercal} \ba_0) \delta(\theta-\sqrt{\lambda_k})$. The whole operation can be approximated by applying an $s$-point FFT on $\{\bx_i, 0\leq i < s\}$, and taking real values. (To be exact, there is a phase factor to be multiplied to the $k$th point in FFT approximation, and is given by $({\sqrt{2 \pi}})^{-1}\eps \exp(-i t_0 k \lambda_{\text{diff}} )$, where we considered the time interval $[t_0, t_0+s \eps]$).

Note that \eqref{eq:LargDiff} is different from the original differential equation in \cite{AvJaSr_Infocom16} where it is $\dot{\bx}(t)= i L \bx(t)$ containing complex coefficients.

\section{Hamiltonian Dynamics and Relation with Quantum Random Walk}
\label{sec:ham_dynamics}
In \cite{AvJaSr_Infocom16} we have studied the Schr\"{o}dinger-type equation of the form \eqref{eq:de_infocom16} with $M$ taken as the adjacency matrix $A$ and $\bpsi$ as the wave function. Now let us consider a similar equation with respect to the graph Laplacian
\be
\label{eq:SrodL}
\dot{\bpsi}(t) = i L \bpsi(t).
\ee
The solution of this dynamics is closely related to the evolution of continuous time quantum random walk and algorithms are developed in \cite{AvJaSr_Infocom16} based on this observation.

Now since the matrix $L$ is real and symmetric, it is sufficient to use the real-imaginary representation of the wave function $\bpsi(t) = \bx(t) + i \by(t)$, $\bx(t), \by(t) \in \mathbb{R}$. Substituting this representation into equation (\ref{eq:SrodL}) and taking real and imaginary parts, we obtain the following system of equations
\bea
\dot{\bx}(t) & = & -L \by(t), \nn \\
\dot{\by}(t) & = & L \bx(t), \nn
\eea
or equivalently in the matrix form
\be
\label{eq:HamDiff}
\frac{d}{dt} \left[\begin{array}{c} \bx(t) \\ \by(t) \end{array}\right] =
\left[\begin{array}{cc} 0 & -L \\ L & 0 \end{array}\right]
\left[\begin{array}{c} \bx(t) \\ \by(t) \end{array}\right].
\ee
Such a system has the following Hamiltonian function
\be
{\cal H} = \frac{1}{2} \bx^\intercal L \bx + \frac{1}{2} \by^\intercal L \by.
\label{eq:ham_schroding}
\ee
Next, the very helpful decomposition
\[
\left[\begin{array}{cc} 0 & -L \\ L & 0 \end{array}\right] =
\left[\begin{array}{cc} 0 & 0 \\ L & 0 \end{array}\right] +
\left[\begin{array}{cc} 0 & -L \\ 0 & 0 \end{array}\right],
\]
together with the observation that
\[
\exp \left( \left[\begin{array}{cc} 0 & 0 \\ L & 0 \end{array}\right] \right) =
\left[\begin{array}{cc} I & 0 \\ L & I \end{array}\right],
\]
leads us to another modification of the leapfrog method known as symplectic split operator algorithm \cite{blanes2006symplectic}: Initialize with
\[
\bdelta \by = -L \bx_0,
\]
then perform the iterations
\bea
\by_{i-1/2} & = & \by_{i-1} - \frac{\eps}{2} \bdelta \by, \label{eq:ham_alg1_0} \\
\bx_i & = & \bx_{i-1} - \eps L \by_{i-1/2}, \label{eq:ham_alg1}
\eea
and update
\bea
\bdelta \by & = & -L \bx_i, \label{eq:ham_alg2_0}\\
\by_i & = & \by_{i-1/2} - \frac{\eps}{2} \bdelta \by. \label{eq:ham_alg2}
\eea
The above modified leapfrog method belongs to the class of symplectic integrator (SI) methods \cite{blanes2008splitting,leimkuhler2004simulating,iserles2009first}. We name the above algorithm as {\it order-2 SI}.

The Hamiltonian system approach can be implemented in two ways:
\begin{enumerate}
\item Form the complex vector $\bx_k+i\by_k$ at each of the $\eps$ intervals. Then $\{\bx_k+i\by_k, 0\leq k < s\}$ with $\bx_0 = \ba_0$ and $\by_0 = \bb_0$ approximates $\exp(iLt)(\ba_0+i\bb_0)$ at $t = 0, \eps,\ldots,(s-1)\eps$ intervals. A direct application of $s$ point FFT with appropriate scaling will give the spectral decomposition as in \eqref{eq:dirac_decompsn}.
\item Note that the formulation in \eqref{eq:HamDiff} is equivalent to the following differential equations
\be
\ddot{\by}(t) + L^2 \by(t) = 0,\quad \ddot{\bx}(t) + L^2 \bx(t) = 0,
\label{eq:Hamiltonian_2nd_tech}
\ee
which are similar to the one in \eqref{eq:LargDiff} except the term $L^2$. Now on the same lines of analysis as in the previous section, taking the real and positive spectrum of just $\by$ component will give $\pi \sum_{k=1}^n \bu_k (\bu_k^{\intercal} \ba_0) \delta(\theta-\lambda_k)$.
\end{enumerate}

\subsection{Fourth order integrator}
The Hamiltonian ${\cal H}$ in \eqref{eq:ham_schroding} associated to the Schr\"odinger-type equation has a special characteristic that it is seperable into two quadratic forms, which help to develop higher order integrators. The $r$ stage integrator has the following form. Between $t$ and $t+\eps$ intervals, we run for $j =1,\ldots,r$,
\bea
\by_j &=& \by_{j-1} + p_j \eps L \bx_{j-1} \nn \\
\bx_j &=& \bx_{j-1} - q_j \eps L \by_j.\nn
\eea
In order to make $q$th order integrator $r\leq q$. For our numerical studies we take the optimized coefficients for order-$4$ derived in \cite{gray1996symplectic}. We call such algorithm as {\it order-4 SI}.

\section{Distributed Implementation}
\label{sec:distributed_impln}

The order-2 symplectic integrator algorithm in \eqref{eq:ham_alg1_0}-\eqref{eq:ham_alg2} can be implemented in a distributed fashion such that each node needs to communicate only to its neighbors. The matrix-vector multiplications in \eqref{eq:ham_alg1} and \eqref{eq:ham_alg2_0}, during one iteration of the algorithm, require diffusion of packets or fluids to the neighbors of every node, and fusion of the received fluids from all the neighbors at each node. Each iteration of the algorithm subsequently has two diffusion-fusion cycles and three synchronization points. A diffusion-fusion cycle consists of $|E|$ packets sent in parallel and hence total number of packets exchanged in one iteration of the algorithm is $2|E|$. Since order-2 SI does not require orthonormalization (unlike classical power iteration and inverse iteration methods for computing eigenelements), and the diffusion-fusion cycle is within one hop neighborhood, time delay of the algorithm will not be too significant. The synchronization points definitely pose some constraints, and demand extra resources.
We have also considered an asynchronous version of the distributed algorithm, which will be presented in
the extended version of the work.

\section{Numerical Results}
\label{sec:numerical_res}
The parameters $\eps$ and $s$ are chosen in the numerical studies satisfying the constraints in \eqref{eq:ft_constraints}. We assume that the maximum degree is known to us.

Note that if the only purpose is to detect eigenvalues, not to compute the eigenvectors, then instead of taking real part of the FFT in the Hamiltonian solution, it is clearly better to compute the absolute value of the complex quantity to get higher peaks. But in the following simulations we look for eigenvectors as well.

For the numerical studies, in order to show the effectiveness of the distributed implementation, we focus on one particular node and plot the spectrum observed at that node. In the plots, $f_{\theta}(k)$ indicates the approximated spectrum at frequency $\theta$ observed on node $k$.

\subsection{Les Mis{\'e}rables network}
In Les Mis{\'e}rables network, nodes are the characters in the well-known novel with the same name and edges are formed if two characters appear in the same chapter. The number of nodes is $77$ and number of edges is $254$. We look for the spectral plot at a specific node called Valjean (with node ID $11$), a character in the associated novel.

\begin{figure}[!htb]
\centering
\subcaptionbox{Euler method\label{fig:lesmiserables_eu_traj}}{\includegraphics[scale=0.525]{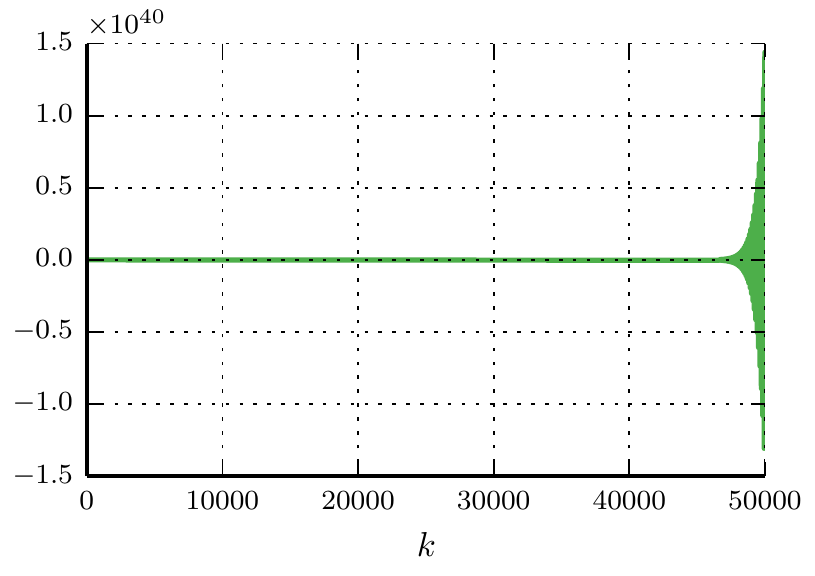}}
\subcaptionbox{Hamiltonian order-2 SI\label{fig:lesmiserables_hm_traj}}{\includegraphics[scale=0.525]{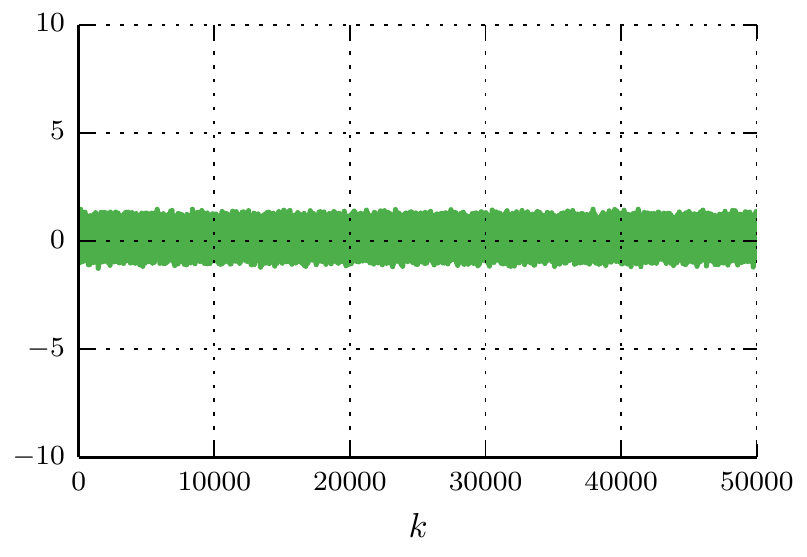}}
\caption{Trajectories}
\label{fig:lesmiserables_trajectories}
\end{figure}

\begin{figure}[h!]
\centering
\includegraphics[scale=0.79]{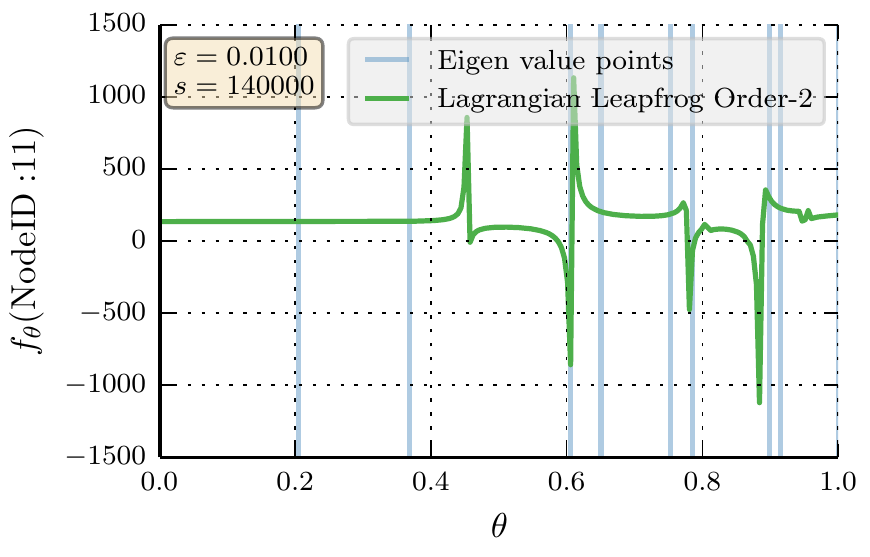}
\vspace{-0.4 cm}
\caption{Les Mis{\'e}rables network: Order-2 Leapfrog}
\label{fig:lesmiserables_or_2_lag}
\end{figure}

\begin{figure}[h!]
\centering
\includegraphics[scale=0.79]{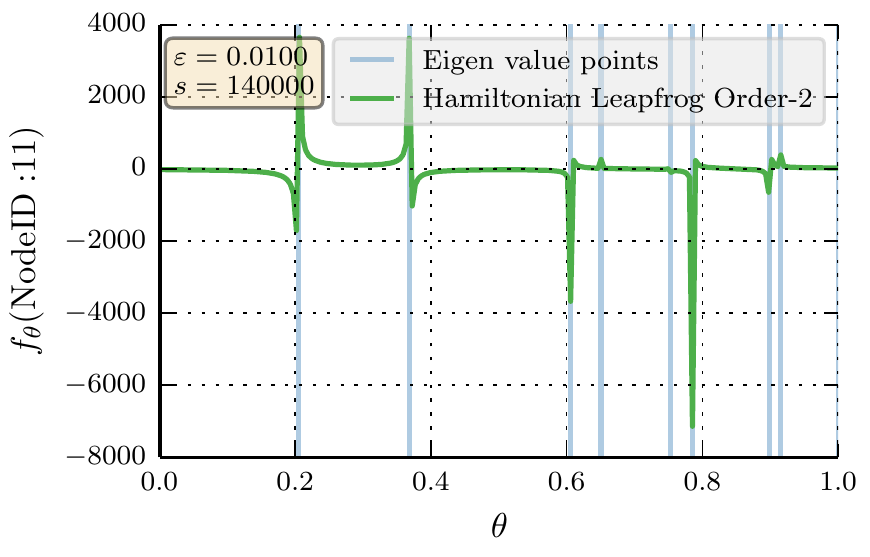}
\vspace{-0.4 cm}
\caption{Les Mis{\'e}rables network: Order-2 SI}
\label{fig:lesmiserables_or_2_ham}
\end{figure}

\begin{figure}[h!]
\centering
\includegraphics[scale=0.79]{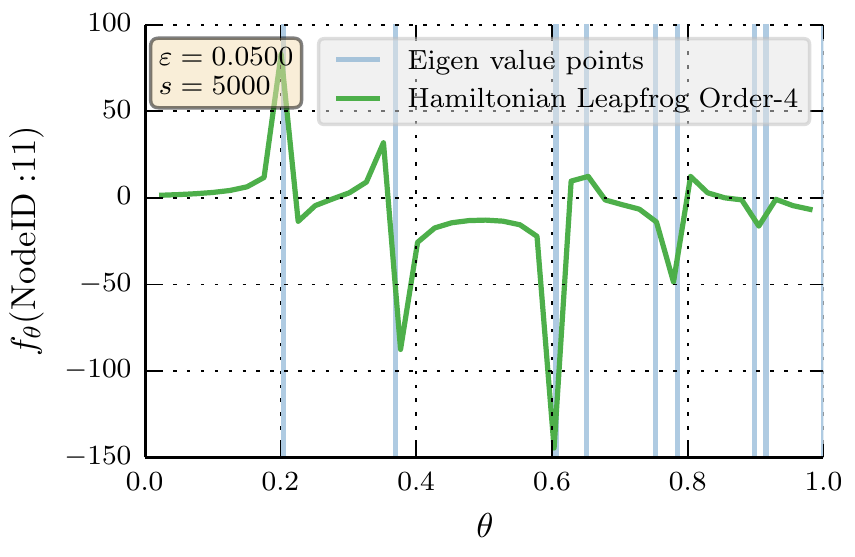}
\vspace{-0.4 cm}
\caption{Les Mis{\'e}rables network: Order-4 SI}
\label{fig:lesmiserables_or_4}
\vspace*{-0.5 cm}
\end{figure}

The instability of the Euler method is clear from Figure~\ref{fig:lesmiserables_eu_traj}, whereas Figure \ref{fig:lesmiserables_hm_traj} shows the guaranteed stability of Hamiltonian SI. Here the y-axis represents the absolute value of $\psi(t)$. (Note the difference in the y-axis scale in the figures). Figure \ref{fig:lesmiserables_or_2_lag} shows the result given by the Lagrangian Leapfrog method from Section~\ref{sec:mechanical_lag}. It can be observed that very few smallest eigenvalues are detected using order-2 Leapfrog compared to the SI technique (order-2) in Figure \ref{fig:lesmiserables_or_2_ham}.
This demonstrates the superiority of the Hamiltonian system approach. Figure \ref{fig:lesmiserables_or_4} shows order-4 SI with much less number of iterations. The precision in order-4 plot can be significantly improved further by increasing the number of iterations.

%
%
%
\vspace*{-0.75 cm}
\vspace*{-0.75 cm}

\subsection{Coauthorship graph in network science}
\begin{figure}[h!]
\hspace{0 cm}
\centering
\includegraphics[scale=0.79]{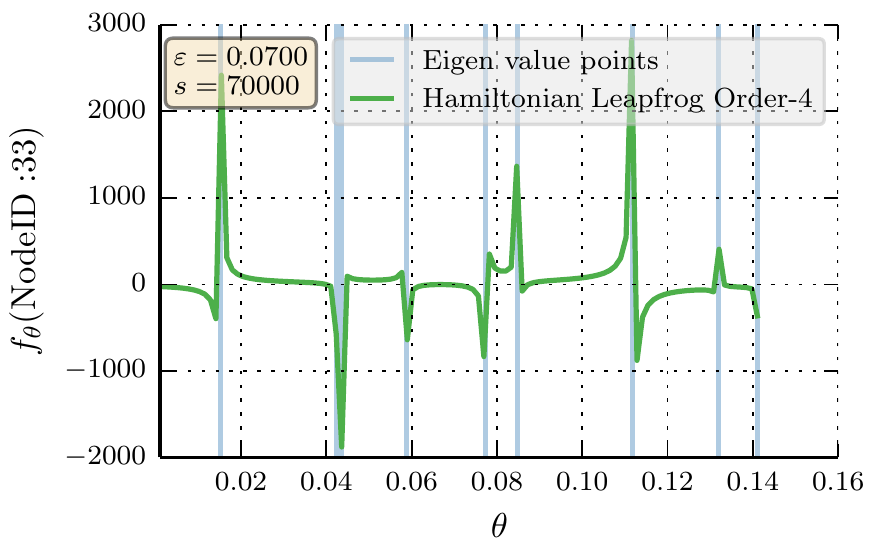}
\vspace{-0.4 cm}
\caption{Coauthorship graph in Network Science: Order-4 SI}
\label{fig:netscience_or_4}
\end{figure}

The coauthorship graph represents a collaborative network of scientists working in network science as compiled by M.\ Newman \cite{Newman_netscience_data}. The numerical experiments are done on the largest connected component with $n=379$ and $m=914$. Figure \ref{fig:netscience_or_4} displays the order-4 SI simulation and it can be seen that even though the eigenvalues are very close, the algorithm is able to distinguish them clearly.
%
\subsection{Enron email network}
\begin{figure}[htb]
\vspace*{-0.4 cm}
\hspace{0 cm}
\centering
\includegraphics[scale=0.79]{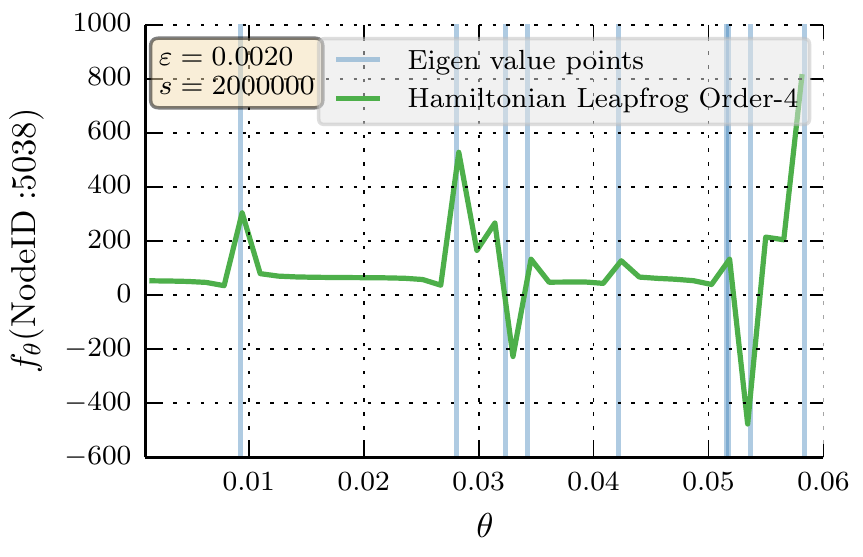}
\vspace{-0.35 cm}
\caption{Enron graph: Order-4 SI}
\label{fig:enron_or_4}
\end{figure}

The nodes in this network are the email ID's of the employees in a company called Enron and the edges are formed when two employees communicated through email\footnote{Data collected from SNAP: \url{http://snap.stanford.edu/data}}. Since the graph is not connected, we take the largest connected component with $33,696$ nodes and $180,811$ edges. The standard MATLAB procedures for eigenelements computation have difficulty to cope with such network sizes. The node under focus is the highest degree node in that component. Simulation result is shown in Figure \ref{fig:enron_or_4}.

\vspace*{-0.4 cm}

\section{Applications}
\label{sec:applications}
In this section we consider two related applications of the distributed
spectrum computation to multi-agent and multi-dimensional (nD) systems.

First, we consider the consensus protocol in the wireless sensor networks \cite{AEN11}.
We model a wireless sensor network as a random geometric graph, with nodes corresponding
to agents (sensors) located on the unit square of the $\mathbb{R}^2$ plane and edges correspond
to possible communication links within radius $R$. The consensus protocol is described
as follows: let $x_k(t)$ be the value of sensor $k$ at time slot $t$,
$$
x_k(t+1) = \sum_{\ell \in N_{[k]}} w_{k \ell} x_{\ell}(t), \quad k=1,...,n,
$$
or in the matrix form $x(t+1)=Wx(t)$, where $N_{[k]}$ is the set of neighbour nodes of node $k$ including the node $k$ itself, and $W=[w_{k \ell}]$ is a matrix of edge weights.
Let the initial value of sensor $k$ be $x_k(0)=m_k$. Then, if the consensus protocol
converges, we have
$$
\lim_{t \to \infty} x(t) = \bar{m}{\bf 1}, \quad \bar{m}=\frac{1}{n}\sum_{k=1}^{n}m_k.
$$
It has been demonstrated in \cite{AEN11} that performance of the best constant consensus protocol  is very
good in sparse wireless networks. The optimal weight value is given by \cite{xiao2004fast}
$$
w_{k \ell} = \frac{2}{\lambda_1(L)+\lambda_{n-1}(L)},
$$
where $\lambda_t(L)$ denotes the $t$-th largest eigenvalue of the graph Laplacian $L=D-W$,
$D=\text{diag}(W{\bf 1})$. Using our distributed approach for the eigenvalues computation, we can
propose a completely distributed, self-tuning, best constant consensus protocol.

We note that in the above example of the consensus protocol, the perfect consensus is
actually reached only in the limit. Often, in practice we would like to obtain consensus
in finite time. The finite-time consensus problem is very challenging and there is no
simple solution for its design (see e.g., \cite{ENA13}). In \cite{Meng12} it has been suggested
to use Iterative Learning Control (ICL) to achieve finite-time consensus. As in \cite{Meng12}
we assume that each agent $k$ can be described by a fairly general Markovian dynamics
$$
y_k(t,r) = g_k(t) + h_k(q)u_k(t,r), \quad k=1,...,n,
$$
where $t$ indicates the time slots, whereas $r$ indicates the ILC iterations.
Here $g_k(t)$ is the zero input responsive function, $h_k(q)$ is the transfer operator with Markovian parameters
and $u_k(t,r)$ is the control input.

The authors of \cite{Meng12} have proposed the following update rule for the control
$$
u_k(t,r+1) = u_k(t,r) + \gamma_k \sum_{\ell \in \CN(k)} a_{k \ell}(y_{\ell}(T,r)-y_k(T,r)),
$$
with $a_{k \ell}$ being elements of the graph adjacency matrix and $\gamma_k$ being the learning gains to be designed,
and shown that under such update rule, the system ``learns'' finite-time consensus, i.e.,
$$
\lim_{r \to \infty} (y_{\ell}(T,r)-y_k(T,r)) = 0, \quad \forall k, \ell \in \{1,...,n\},
$$
if the following condition holds
$$
\rho(I-H\Gamma L) < 1,
$$
where $\Gamma=\text{diag}\{\gamma_1,...,\gamma_n\}$ is the diagonal matrix of gains,
$H=\text{diag}\{h_1(T),...,h_n(T)\}$ is the diagonal matrix of the response function at time $T$
and $\rho(A)$ is the spectral radius of matrix $A$. In fact, if the communication graph
is undirected and connected, the condition $\rho(I-H\Gamma L) < 1$ is always satisfied.
However, in practice, it is good to be not too aggressive in learning \cite{Longman2002},
and hence our distributed procedure can be used to choose gains $\{\gamma_k\}$ in such a way so that
the value $\rho(I-H\Gamma L)$ estimated by our algorithm will be not too small and not
too close to one.

\section{Conclusions and future research}
\label{sec:conclusions}
We have proposed a distributed approach for the eigenvalue-eigenvector problem for
graph matrices based on Hamiltonian dynamics, symplectic integrators and smoothed Fourier
transform. We demonstrate with various network sizes that the proposed approach efficiently
scales and finds, with higher resolution, closely situated eigenvalues and associated
eigenvectors of graph matrices.

In the future we hope to present asynchronous versions of the introduced algorithms,
to extend the proposed approaches to some classes of non-symmetric matrices
and to design an automatic or a semi-automatic procedure
for the identification of dominant eigenvalues and eigenvectors.


\section*{Acknowledgement}

This work is supported by INRIA Bell Labs joint lab (ADR Network Science). We also would like to thank
Leonid Freidovich for stimulating discussions. This is the author version of the IEEE nDS 2017 article.

\bibliographystyle{IEEEtran}
\bibliography{IEEEabrv,nds17_v1}
\end{document}